\documentclass{article}
\usepackage{amsfonts}
\usepackage{amsmath}

\newtheorem{theorem}{Theorem}

\newtheorem{corollary}[theorem]{Corollary}

\newtheorem{lemma}[theorem]{Lemma}

\begin{document}

\title{How many families survive for a long time\thanks{%
This work is supported by the RSF under a grant 14-50-00005. } }
\author{Vatutin V.A.\thanks{%
Department of Discrete Mathematics, Steklov Mathematical Institute, 8,
Gubkin str., 119991, Moscow, Russia; e-mail: vatutin@mi.ras.ru}, Dyakonova
E.E.\thanks{%
Department of Discrete Mathematics, Steklov Mathematical Institute, 8,
Gubkin str., 119991, Moscow, Russia; e-mail: elena@mi.ras.ru}}
\date{}
\maketitle

\begin{abstract}
A critical branching process $\left\{ Z_{k},k=0,1,2,...\right\} $ in a
random environment generated by a sequence of independent and identically
distributed random reproduction laws is considered.\ Let $Z_{p,n}$ be the
number of particles at time $p\leq n$ having a positive offspring number at
time $n$. \ A theorem is proved describing the limiting behavior, as $%
n\rightarrow \infty $ of the distribution of a properly scaled process $\log
Z_{p,n}$ under the assumptions $Z_{n}>0$ and $p\ll n$.
\end{abstract}

\medskip


\medskip 


\section{Introduction}

We consider branching processes in random environment specified by sequences
of independent identically distributed random laws. Denote by $\Delta $ the
space of proper probability measures on $\mathbb{N}_{0}=\left\{
0,1,2,...\right\} ={0}\cup \mathbb{N}_{+}$. Let $Q$ be a random variable
taking values in $\Delta $. An infinite sequence
\begin{equation}
\Pi =(Q_{1},Q_{2},\ldots )  \label{DefEnvir}
\end{equation}%
of i.i.d. copies of $Q$ is said to form a \emph{random environment}. In the
sequel we make no difference between the laws
\begin{equation*}
Q=(Q(\{0\}),Q(\{1\}),...,Q(\{k\}),...),\,Q_{n}=(Q_{n}(\{0\}),Q_{n}(\{1%
\}),...,Q_{n}(\{k\}),...)
\end{equation*}%
and the generating functions%
\begin{equation}
f\left( s\right) =f(s;Q):=\sum_{k=0}^{\infty }Q(\{k\})s^{k},\quad
f_{n}\left( s\right) =f_{n}(s;Q):=\sum_{k=0}^{\infty
}Q_{n}(\{k\})s^{k},\quad n\in \mathbb{N}_{+}.  \label{genfuu}
\end{equation}%
A sequence of $\mathbb{N}_{0}$-valued random variables $\mathcal{Z}:\mathcal{%
=}\left( Z_{0},Z_{1},\ldots \right) $ is called a \emph{branching process in
the random environment} $\Pi $, if $Z_{0}$ is independent of $\Pi $ and
\begin{equation*}
\mathbf{E}\left[ s_{n}^{Z_{n}}\;|\;Z_{0},Z_{1},...,Z_{n-1},Q_{1},Q_{2},%
\ldots ,Q_{n}\right] \ =\left( f_{n}\left( s\right) \right)
^{Z_{n-1}},\;n\in \mathbb{N}_{+}.
\end{equation*}

Thus, $Z_{n}$ is the $n$th generation size of the population and $Q_{n}$ is
the distribution of the number of children of an individual at generation $%
n-1$. We assume that $Z_{0}=1$ for convenience and denote the corresponding
probability measure on the underlying probability space by~$\mathbf{P}$. (If
we refer to other probability spaces, then we use notation $\mathbb{P}$, $%
\mathbb{E}$ \ and $\mathbb{L}$ for the respective probability measures,
expectations and laws.)

As it turns out the properties of $\mathcal{Z}$ are mainly determined by its
associated random walk $\mathcal{S}:=\left\{ S_{n},n\geq 0\right\} $. This
random walk has initial state $S_{0}=0$ and increments $X_{n}=S_{n}-S_{n-1}$%
, $n\geq 1,$ defined as
\begin{equation*}
X_{n}\ :=\log f_{n}^{\prime }(1),
\end{equation*}%
which are i.i.d. copies of the logarithmic mean offspring number
\begin{equation*}
X\ :=\ \log f^{\prime }(1).
\end{equation*}

Following \cite{AGKV05} we call the process $\mathcal{Z}$ \textit{critical}
if and only if the random walk $\mathcal{S}$ is oscillating, that is,
\begin{equation*}
\limsup_{n\rightarrow \infty }S_{n}=\infty \ \text{ and }\
\liminf_{n\rightarrow \infty }S_{n}=-\infty
\end{equation*}%
with probability 1.

It is shown in \cite{AGKV05} that the extinction moment of the critical
branching process in a random environment is finite with probability $1$.
Moreover, it is established in \cite{AGKV05} that if
\begin{equation}
\lim_{n\rightarrow \infty }\mathbf{P}\left( S_{n}>0\right) =:\rho \in (0,1),
\label{Spit}
\end{equation}%
then, as $n\rightarrow \infty $ (and some mild additional assumptions to be
specified later on)
\begin{equation}
\mathbf{P}(Z_{n}>0)\sim \theta \mathbf{P}(\min \left(
S_{0},S_{1},...,S_{n}\right) \geq 0)=\theta n^{\rho -1}l(n),  \label{tails}
\end{equation}%
where $l(n)$ is a slowly varying function and $\theta $ is a known positive
constant whose explicit expression is given by formula (\ref{DefTheta})
below.

Let
\begin{equation*}
\mathcal{A}:=\{0<\alpha <1;\,\left\vert \beta \right\vert <1\}\cup
\{1<\alpha <2;\left\vert \beta \right\vert <1\}\cup \{\alpha =1,\beta
=0\}\cup \{\alpha =2,\beta =0\}
\end{equation*}%
be a subset in $\mathbb{R}^{2}.$ For $(\alpha ,\beta )\in \mathcal{A}$ and a
random variable $X$ write $X\in \mathcal{D}\left( \alpha ,\beta \right) $ if
the distribution of $X$ belongs to the domain of attraction of a stable law
with characteristic function%
\begin{equation}
\mathcal{H}_{\alpha ,\beta }\mathbb{(}t\mathbb{)}:=\exp \left\{
-c|t|^{\,\alpha }\left( 1+i\beta \frac{t}{|t|}\tan \frac{\pi \alpha }{2}%
\right) \right\} ,\ c>0,  \label{std}
\end{equation}%
and, in addition, $\mathbf{E}X=0$ if this moment exists.

Denote by $\left\{ c_{n},n\in \mathbb{N}_{+}\right\} $ a sequence of
positive integers specified by the relation%
\begin{equation*}
c_{n}:=\inf \left\{ u\geq 0:G(u)\leq n^{-1}\right\} ,
\end{equation*}%
where
\begin{equation*}
G(u):=\frac{1}{u^{2}}\int_{-u}^{u}x^{2}\mathbf{P}(X\in dx).
\end{equation*}%
It is known (see, for instance, \cite[Ch. XVII, \S 5]{FE}) that, for every $%
X\in \mathcal{D}(\alpha ,\beta )$ the function $G(u)$ is regularly varying
with index $(-\alpha )$, and, therefore,
\begin{equation}
c_{n}=n^{1/\alpha }l_{1}(n),\ n\in \mathbb{N}_{+},  \label{asyma}
\end{equation}%
for some function $l_{1}(n),$ slowly varying at infinity. In addition, the
scaled sequence $\left\{ S_{n}/c_{n},\,n\geq 1\right\} $ converges in
distribution, as $n~\rightarrow ~\infty ,$ to the stable law given by (\ref%
{std}). \

Observe that if $X\in \mathcal{D}\left( \alpha ,\beta \right) ,$ then the
quantity $\rho $ in~(\ref{Spit}) is calculated by the formula (see,$\ $for
instance, \cite{Zol57})
\begin{equation}
\displaystyle\rho =\left\{
\begin{array}{ll}
\frac{1}{2},\ \text{if \ }\alpha =1\text{ or }2, &  \\
\vspace{0.1cm} &  \\
\frac{1}{2}+\frac{1}{\pi \alpha }\arctan \left( \beta \tan \frac{\pi \alpha
}{2}\right) ,\text{ otherwise}. &
\end{array}%
\right.  \label{ro}
\end{equation}%
In particular, $\rho \in \left( 0,1\right) $.

Denote%
\begin{equation*}
M_{n}:=\max \left( S_{1},...,S_{n}\right) ,\quad L_{k,n}:=\min_{k\leq j\leq
n}S_{j},\quad L_{n}:=L_{0,n}=\min \left( S_{0},S_{1},...,S_{n}\right)
\end{equation*}%
and introduce two functions
\begin{eqnarray*}
V(x) &:&=1+\sum_{k=1}^{\infty }\mathbf{P}\left( -S_{k}\leq x,M_{k}<0\right)
,\ x\geq 0, \\
U(x) &:&=1+\sum_{k=1}^{\infty }\mathbf{P}\left( -S_{k}>x,L_{k}\geq 0\right)
,\ x\leq 0,
\end{eqnarray*}%
and $0$ elsewhere. In particular, $V(0)=U(0)=1$.

The fundamental properties of $V,U$ are the identities
\begin{equation}
\begin{array}{rl}
\mathbf{E}[V(x+X);x+X\geq 0]\ =\ V(x)\ , & x\geq 0\ , \\
\mathbf{E}[U(x+X);x+X<0]\ =\ U(x)\ , & x\leq 0\ ,%
\end{array}
\label{harm}
\end{equation}%
which hold for any oscillating random walk.

It follows from (\ref{harm}) that $V$ and $U$ give rise to probability
measures $\mathbf{P}_{x}^{+}$ and $\mathbf{P}_{x}^{-}$ being important for
subsequent arguments. The construction procedures of these probability
measures are explained for $\mathbf{P}_{x}^{+}$ in \cite{AGKV05} and for $%
\mathbf{P}_{x}^{-}$ in \cite{ABKV2009} and \cite{VD2003} in detail. We only
recall that if $T_{0},T_{1},\ldots $ is a sequence of random variables and
the random walk $\mathcal{S}=(S_{n},n\geq 0)$ with $S_{0}=x$ are both
adapted to some filtration $\mathfrak{F}=(\mathcal{F}_{n},n\geq 0)$, then,
for each fixed $n$ and a bounded and measurable function $g_{n}:\mathbb{R}%
^{n+1}\rightarrow \mathbb{R}$ the measures above are specified by the
equalities%
\begin{eqnarray*}
\mathbf{E}_{x}^{+}[g_{n}(T_{0},\ldots ,T_{n})]\ &:&=\ \frac{1}{V(x)}\mathbf{E%
}_{x}[g_{n}(T_{0},\ldots ,T_{n})V(S_{n});L_{n}\geq 0], \\
\mathbf{E}_{x}^{-}[g_{n}(T_{0},\ldots ,T_{n})]\ &:&=\ \frac{1}{U(x)}\mathbf{E%
}_{x}[g_{n}(T_{0},\ldots ,T_{n})U(S_{n});M_{n}<0].
\end{eqnarray*}%
Observe that under the measure $\mathbf{P}_{x}^{+}$ the sequence $%
S_{0},S_{1},\ldots $ is a Markov process with state space $\mathbb{R}^{0}:=[0,\infty )$ and
transition probabilities
\begin{equation*}
\mathbf{P}^{{+}}(x,dy)\ :=\ \frac{1}{V(x)}\mathbf{P}\{x+X\in dy\}V(y)\
,\quad x\geq 0,y\geq 0.
\end{equation*}%
It is the random walk conditioned never to enter $\mathbb{R}^{-}:=(-\infty ,0),$ while under
the measure $\mathbf{P}_{x}^{-}$ the process $S_{0},S_{1},\ldots $ becomes a
Markov chain with state space $\mathbb{R}^{-}$ and transition kernel
\begin{equation}
\mathbf{P}^{{-}}(x,dy)\ :=\ \frac{1}{U(x)}\mathbf{P}\{x+X\in dy\}U(y),\quad
x\leq 0,y<0.  \label{DefPminus}
\end{equation}%
Note that $\mathbf{P}^{{-}}(x,[0,\infty ))=0$, thus the Markov process never
enters $[0,\infty )$ again. It may, however, start from the boundary $x=0$.

We now describe in brief a construction of L\'evy processes conditioned to
stay nonnegative following basically the definitions given in \cite{Chau06}
and \cite{CCh2006}.

Let $\Omega :=D\left( [0,\infty ),\mathcal{R}\right) $ be the space of
real-valued c\`{a}dl\`{a}g paths on the real half-line $[0,\infty )$ and let
$\mathcal{B}:=\left\{ B_{t},t\geq 0\right\} $ be the coordinate process
defined by the equality $B_{t}(\omega )=\omega _{t}$ for $\omega \in \Omega $%
. In the sequel we consider also the spaces $\Omega _{U}:=D\left( [0,U],%
\mathcal{R}\right) ,U>0.$

We endow the spaces $\Omega $ and $\Omega _{U}$ with Skorokhod topology and
denote by $\mathcal{F}=\left\{ \mathcal{F}_{t},t\geq 0\right\} $ and by $%
\mathcal{F}^{U}\mathcal{=}\left\{ \mathcal{F}_{t},t\in \lbrack 0,U]\right\} $
(with some misuse of notation) the natural filtrations of the processes $%
\mathcal{B}$ and $\mathcal{B}^{U}=\left\{ B_{t},t\in \lbrack 0,U]\right\} $.

Let $\mathbb{P}_{x}$ be the law on $\Omega $ of an $\alpha -$stable L\'evy
process $\mathcal{B}$, $\alpha \in (0,2]$ started at $x$ and let $\mathbb{P}=%
\mathbb{P}_{0}$. Denote by $\rho :=\mathbb{P}\left( B_{1}\geq 0\right) $ the
positivity parameter of the process $\mathcal{B}$ (in fact, this quantity is
the same as in (\ref{ro}) for a random walk whose increments $X_{i}\in
D(\alpha ,\beta )$). We now introduce an analogue of the measure $\mathbf{P}%
^{+}$ for L\'evy processes. Namely, following \cite{Chau97} we specify for all
$t>0$ and events $\mathcal{A}\in \mathcal{F}_{t}$ the law $\mathbb{P}_{x}^{+}$ on $%
\Omega $ \ of the L\'evy process starting at point $x>0$ and conditioned to
stay nonnegative by the equality
\begin{equation*}
\mathbb{P}_{x}^{+}\left( \mathcal{A}\right) :=\frac{1}{x^{\alpha \left(
1-\rho \right) }}\mathbb{E}_{x}\left[ B_{t}^{\alpha \left( 1-\rho \right)
}I\left\{ \mathcal{A}\right\} I\left\{ \inf_{0\leq u\leq t}B_{u}\geq
0\right\} \right] ,
\end{equation*}%
where $I\left\{ \mathcal{C}\right\} $ is the indicator of the event $%
\mathcal{C}$.

Thus, $\mathbb{P}_{x}^{+}$ is an $h-$transform of the L\'evy process killed
when it first enters the negative half-line. The corresponding positive
invariant function is $H(x)=x^{\alpha \left( 1-\rho \right) }$ in this case.

This definition has no sense for $x=0$. However, it is shown in \cite{Chau06}
that one can construct a law $\mathbb{P}^{+}=\mathbb{P}_{0}^{+}$ and a c\`{a}%
dl\`{a}g Markov process with the same semigroup as $\left( \mathcal{B}%
,\left\{ \mathbb{P}_{x}^{+},x>0\right\} \right) $ in such a way that $%
\mathbb{P}^{+}\left( B_{0}=0\right) =1$. Moreover,%
\begin{equation*}
\mathbb{P}_{x}^{+}\Longrightarrow \mathbb{P}^{+},\text{ as }x\downarrow 0%
\text{.}
\end{equation*}%
Here and in what follows the symbol $\Longrightarrow $ stands for the weak
convergence in the respective space of c\`{a}dl\`{a}g functions endowed with
the Skorokhod topology.

\bigskip Let $\mathbb{P}^{(m)}$ be the law on $\Omega _{1}$ of the meander
of length 1 associated with $\left( \mathcal{B},\mathbb{P}\right) ,$ i.e.
\begin{equation}
\mathbb{P}^{(m)}\left( \cdot \right) :=\lim_{x\downarrow 0}\mathbb{P}%
_{x}\left( \cdot \Big|\inf_{0\leq u\leq 1}B_{u}\geq 0\right) .  \label{MesPm}
\end{equation}%
Thus, the law $\mathbb{P}^{(m)}$ may be viewed as the law of the L\'evy
process $\left( \mathcal{B},\mathbb{P}\right) $ conditioned to stay
nonnegative on the time-interval $(0,1]$ while the law $\mathbb{P}^{+}$
specified earlier corresponds to the law of the L\'evy process conditioned to
stay nonnegative on the whole time interval $(0,\infty )$.

As shown in \cite{Chau06}, $\mathbb{P}^{(m)}$ and $\mathbb{P}^{+}$ are
absolutely continuous with respect to each other: for every event $\mathcal{A%
}\in \mathcal{F}_{1}$
\begin{equation}
\mathbb{P}^{+}\left( \mathcal{A}\right) =C_{0}\mathbb{E}^{(m)}\left[
I\left\{ \mathcal{A}\right\} B_{1}^{\alpha \left( 1-\rho \right) }\right] ,
\label{DefPplus}
\end{equation}%
where (see, for instance, formulas (3.5)-(3.6) and (3.11) in \cite{CCh2006})
\begin{equation}
C_{0}:=\lim_{n\rightarrow \infty }V(c_{n})\mathbf{P}\left( L_{n}\geq
0\right) =\left( \mathbb{E}^{(m)}\left[ B_{1}^{\alpha \left( 1-\rho \right) }%
\right] \right) ^{-1}\in \left( 0,\infty \right) ,  \label{AsH}
\end{equation}%
and $\mathbb{E}^{(m)}$ is the expectation with respect to $\mathbb{P}^{(m)}.$
In fact, one may extend the absolute continuity given in (\ref{DefPplus}) to
an arbitrary interval $[0,U]$ be considering the respective space $\Omega
_{U}$ instead of $\Omega _{1}$ and conditioning by the event $\inf_{0\leq
u\leq U}B_{u}\geq 0$ in (\ref{MesPm}).

We now come back to the branching processes in random environment and set%
\begin{equation*}
\zeta (a):=\frac{\sum_{k=a}^{\infty }k^{2}Q(\{k\})}{\left( f^{\prime }\left(
1\right) \right) ^{2}},\ a\in \mathbb{N}_{0}.
\end{equation*}

In what follows we say that

1) \textit{Condition }$A1$\textit{\ is valid} if $X\in \mathcal{D}\left(
\alpha ,\beta \right) ;$

2) \textit{Condition }$A2$\textit{\ is valid} if
\begin{equation*}
\mathbf{E}\left( \log ^{+}\zeta (a)\right) ^{\alpha +\varepsilon }<\infty
\end{equation*}%
for some $\varepsilon >0$ and $a\in \mathbb{N}_{0}$;

3) \textit{Condition }$A$\textit{\ is valid} if Conditions $A1$- $A2$ hold
true and, in addition, the parameter $p=p(n)$ tends to infinity as $%
n\rightarrow \infty $ in such a way that%
\begin{equation}
\lim_{n\rightarrow \infty }n^{-1}p=\lim_{n\rightarrow \infty }n^{-1}p(n)=0.
\label{Msmalln}
\end{equation}

Branching processes in random environment meeting Condition $A$ have been
investigated in a recent paper \cite{VD2016}. The paper includes a
Yaglom-type functional limit theorem describing the asymptotic properties of
the process $\left\{ Z_{\left[ tp\right] },0\leq t<\infty \right\} $ given $%
Z_{n}>0$.

Here we investigate the structure of the process $\mathcal{Z}$ in more
detail and prove a conditional limit theorem for the so-called reduced
process $\left\{ Z_{p,n},0\leq p\leq n\right\} ,$ where $Z_{p,n}$ is the
number of particles in the process at time $p\in \lbrack 0,n]$, each of
which has a nonempty offspring at time $n$. Our main result looks as follows.

\begin{theorem}
\label{T_reduced}If Condition $A$ is valid, then for any $x\geq 0$
\begin{eqnarray}
\lim_{n\rightarrow \infty }\mathbf{P}\left( \frac{\log Z_{p,n}}{c_{p}}\geq
x|Z_{n}>0,Z_{0}=1\right) &=&\mathbb{E}^{+}\left[ \left( 1-\frac{x}{B_{1}}%
\right) ^{\alpha \left( 1-\rho \right) }I\left\{ B_{1}\geq x\right\} \right]
\notag \\
&=&\mathbb{P}^{+}\left( \inf_{1\leq v<\infty }B_{v}\geq x\right) .
\label{One_dim}
\end{eqnarray}
\end{theorem}

Here and in what follows $\mathbb{E}^{+}$ is the expectation with respect to
$\mathbb{P}^{+}$\ . This result complements papers \cite{BV} and \cite%
{Vat2002} where it was established (under the assumptions $\mathbf{E}\left[ X%
\right] =0,$ $\sigma ^{2}=\mathbf{E}\left[ X^{2}\right] \in \left( 0,\infty
\right) $ and some additional technical conditions) that, as $n\rightarrow
\infty $
\begin{eqnarray*}
&&\mathcal{L}\left( \left\{ \frac{\log Z_{\left[ nt\right] ,n}}{\sigma \sqrt{%
n}},0\leq t\leq 1\right\} \Big|Z_{n}>0,Z_{0}=1\right) \\
&&\qquad \qquad \qquad \qquad \Longrightarrow \mathbb{L}^{(m)}\left(
\inf_{t\leq v\leq 1}B_{v},0\leq t\leq 1\right) =\mathcal{L}\left(
\inf_{t\leq v\leq 1}B_{v}^{+}\right) ,
\end{eqnarray*}%
where $B^{+}:=\{B_{v}^{+}, 0\leq \leq v\leq 1\}$ is the standard Brownian meander.

The study of the reduced processes has a rather long history. Reduced
processes for ordinary Galton-Watson branching processes where introduced by
Fleischmann and Prehn \cite{FP}, who discussed the subcritical case.
Critical and supercritical reduced Galton-Watson processes have been
investigated by Zubkov \cite{Zub} and Fleischmann and Siegmund-Schultze \cite%
{FZ}. The first results for reduced branching processes in random
environment were obtained by Borovkov and Vatutin \cite{BV} and Fleischmann
and Vatutin \cite{FV99}. They considered (under the annealed approach) the
case when the support of the measure $\mathbf{P}$ is concentrated only on
the set of fractional-linear generating functions. Vatutin proved in \cite%
{Vat2002} a limit theorem for the critical reduced processes under the
annealed approach and general reproduction laws of particles. Papers \cite%
{VD}, \cite{VD2005},\cite{VD2007} and \cite{VD2009} consider the properties of
the critical reduced branching processes in random environment under the
quenched approach (see also survey \cite{VDS2013}).

\section{Auxiliary results}

To prove the main results of the paper we need to know the asymptotic
behavior of the function $V(x)$ as $x\rightarrow \infty $. The following
lemma gives the desired representation.

\begin{lemma}
\label{Renew2} (compare with Lemma 13 in \cite{VW09} and Corollary 8 in \cite%
{Don12}) If $X\in \mathcal{D}\left( \alpha ,\beta \right) ,$ then there
exists a function $l_{0}(x)$ slowly varying at infinity such that
\begin{equation*}
V(x)\sim x^{\alpha (1-\rho )}l_{0}(x)
\end{equation*}%
as $x\rightarrow \infty $.
\end{lemma}

In the sequel we use the symbols $K,K_{1},K_{2},...$ to denote different
constants. They are not necessarily the same in different formulas.

Our next result is a combination (with a slight reformulation) of Lemma 2.1
in \cite{AGKV05} and Corollaries 3 and 8 in \cite{Don12}:

\begin{lemma}
\label{L_minimBelow} If $X\in \mathcal{D}\left( \alpha ,\beta \right) ,$
then (compare with (\ref{tails})), as\ $n\rightarrow \infty $
\begin{equation}
\mathbf{P}\left( L_{n}\geq -r\right) \sim V(r)\mathbf{P}\left( L_{n}\geq
0\right) \sim V(r)n^{\rho -1}l(n)  \label{UnifPrecise}
\end{equation}%
uniformly in $0\leq r\ll c_{n},$ and there exists a constants $K_{1}>0$ such
that%
\begin{equation}
\mathbf{P}\left( L_{n}\geq -r\right) \leq K_{1}V(r)\mathbf{P}\left(
L_{n}\geq 0\right) ,\ r\geq 0,\, n\geq 1.  \label{IntermedEst}
\end{equation}
\end{lemma}

For $U>0$ and $pU\leq n$ let
\begin{eqnarray*}
\mathcal{Q}_{U}^{p,n} &:&=\left\{ \frac{S_{\left[ pu\right] }}{c_{p}},0\leq
u\leq U\right\} ,\quad \mathcal{Q}^{p,n}:=\mathcal{Q}_{\infty }^{p,n}, \\
\mathcal{S}_{U}^{p,n} &:&=\left\{ \frac{S_{[pU]+\left[ \left( n-pU\right) t%
\right] }}{c_{n}},0\leq t\leq 1\right\} ,\quad \mathcal{S}^{n}:=\mathcal{S}%
_{0}^{p,n}.
\end{eqnarray*}

Let $\phi _{1}:\Omega _{1}\rightarrow \mathbb{R}$ be a bounded uniformly
continuous functional and $\left\{ \varepsilon _{n},n\in \mathbb{N}%
_+\right\} $ be a sequence of positive numbers vanishing as $n\rightarrow
\infty $.

\begin{lemma}
\label{L_uniform}(see \cite{VD2016})\ If Condition $A1$ is valid then
\begin{equation*}
\mathbf{E}\left[ \phi _{1}\left( \mathcal{S}^{n}\right) |L_{n}\geq -x\right]
\rightarrow \mathbb{E}^{(m)}\left[ \phi _{1}(\mathcal{B}^{1})\right]
\end{equation*}%
as $n\rightarrow \infty $ uniformly in $0\leq x\leq \varepsilon _{n}c_{n}.$
\end{lemma}

\begin{lemma}
\label{C_funct}(see \cite{VD2016}) If Conditions $A1$ and (\ref{Msmalln})
are valid, then for any $r\geq 0$
\begin{equation*}
\mathcal{L}\left( \mathcal{Q}^{p,n}\Big|L_{n}\geq -r\right) \Longrightarrow
\mathbb{L}^{+}\left( \mathcal{B}\right)
\end{equation*}%
as $n\rightarrow \infty $.
\end{lemma}

For $x\geq 0$ we set for brevity%
\begin{eqnarray}
D(x) &:&=C_{0}\mathbb{E}^{(m)}\left[ \left( B_{1}-x\right) ^{\alpha \left(
1-\rho \right) }I\left\{ B_{1}\geq x\right\} \right]  \notag \\
&=&\mathbb{E}^{+}\left[ \left( 1-\frac{x}{B_{1}}\right) ^{\alpha \left(
1-\rho \right) }I\left\{ B_{1}\geq x\right\} \right] .  \label{defD}
\end{eqnarray}

\begin{lemma}
\label{L_explicit}For any $x>0$%
\begin{equation*}
\mathbb{P}^{+}\left( \inf_{1\leq u<\infty }B_{u}\geq x\right) =D(x).
\end{equation*}
\end{lemma}

\textbf{Proof}. According to Theorem 5 in \cite{Chau97} for any pair $%
0<x\leq y$
\begin{equation*}
\mathbb{P}_{y}^{+}\left( \inf_{0\leq u<\infty }B_{u}\geq x\right) =\left( 1-%
\frac{x}{y}\right) ^{\alpha \left( 1-\rho \right) }I\left\{ x\leq y\right\} .
\end{equation*}%
Hence,%
\begin{eqnarray*}
\mathbb{P}^{+}\left( \inf_{1\leq u<\infty }B_{u}\geq x\right)
&=&\int_{x}^{\infty }\mathbb{P}^{+}\left( B_{1}\in dy\right) \mathbb{P}%
_{y}^{+}\left( \inf_{0\leq u<\infty }B_{u}\geq x\right) \\
&=&\int_{x}^{\infty }\mathbb{P}^{+}\left( B_{1}\in dy\right) \left( 1-\frac{x%
}{y}\right) ^{\alpha \left( 1-\rho \right) }=D\mathbb{(}x\mathbb{)},
\end{eqnarray*}%
as desired.

In the sequel we agree to write $a_{n}\ll b_{n}~$if $\lim_{n\rightarrow
\infty }a_{n}/b_{n}=0.$ In prticular, $\lim_{n\gg p\rightarrow \infty }$
means that the limit of the respective expression is claculated as $%
p,n\rightarrow \infty $ in such a way that $pn^{-1}\rightarrow 0$.

The following statement will be useful for the subsequent proofs.

\begin{lemma}
\label{L_MinCond}If Condition $A1$ is valid then, for any $x\geq 0$
\begin{equation}
\lim_{n\gg p\rightarrow \infty }\mathbf{P}\left( L_{p,n}\geq
xc_{p}|L_{n}\geq -r\right) =D(x).  \label{MinCon}
\end{equation}
\end{lemma}

\textbf{Proof.} We select $N>x$ and write%
\begin{eqnarray}
0 &\leq &\mathbf{P}\left( L_{p,n}\geq xc_{p};L_{n}\geq -r\right) -\mathbf{P}%
\left( S_{p}\leq Nc_{p},L_{p,n}\geq xc_{p};L_{n}\geq -r\right)  \notag \\
&\leq &\mathbf{P}\left( S_{p}>Nc_{p};L_{n}\geq -r\right) .  \label{Decomp}
\end{eqnarray}%
It follows from Lemma \ref{C_funct} that for any $\varepsilon >0$ one can
find $N_{0}=N_{0}\left( \varepsilon \right) $ such that for all $N\geq N_{0}$
\begin{equation}
\mathbf{P}\left( S_{p}>Nc_{p};L_{n}\geq -r\right) \leq \varepsilon \mathbf{P}%
\left( L_{n}\geq -r\right) .  \label{SSmall}
\end{equation}

To proceed further we denote by $\mathcal{S}^{\ast }:=\left( S_{0}^{\ast
}=0,S_{1}^{\ast },...,S_{n}^{\ast },...\right) $ an independent
probabilistic copy of the random walk $\mathcal{S}$ and let
\begin{equation*}
L_{k}^{\ast }:=\min \left( S_{0}^{\ast },S_{1}^{\ast },...,S_{k}^{\ast
}\right) \text{.}
\end{equation*}%
Then, for $N>x\geq 0$%
\begin{eqnarray*}
&&\mathbf{P}\left( S_{p}\leq Nc_{p},L_{p,n}\geq xc_{p};L_{n}\geq -r\right) \\
&&\qquad =\int_{xc_{p}}^{Nc_{p}}\mathbf{P}\left( S_{p}\in dy;L_{p}\geq
-r\right) \mathbf{P}\left( L_{n-p}^{\ast }\geq xc_{p}-y\right) \\
&&\qquad =\int_{x}^{N}\mathbf{P}\left( S_{p}\in c_{p}dz;L_{p}\geq -r\right)
\mathbf{P}\left( L_{n-p}^{\ast }\geq \left( x-z\right) c_{p}\right) .
\end{eqnarray*}%
Since $p\ll n,$ we deduce by (\ref{UnifPrecise}), (\ref{AsH}) and properties
of regularly varying functions that if $n\rightarrow \infty $ then for any $%
\varepsilon >0$
\begin{eqnarray*}
\mathbf{P}\left( L_{n-p}^{\ast }\geq \left( x-z\right) c_{p}\right) &\sim
&V(\left( z-x\right) c_{p})\mathbf{P}\left( L_{n}\geq 0\right) \\
&\sim &\left( z-x\right) ^{\alpha \left( 1-\rho \right) }V(c_{p})\mathbf{P}%
\left( L_{n}\geq 0\right) \\
&=&\left( z-x\right) ^{\alpha \left( 1-\rho \right) }V(c_{p})\mathbf{P}%
\left( L_{p}\geq 0\right) \frac{\mathbf{P}\left( L_{n}\geq 0\right) }{%
\mathbf{P}\left( L_{p}\geq 0\right) } \\
&\sim &\left( z-x\right) ^{\alpha \left( 1-\rho \right) }C_{0}\frac{\mathbf{P%
}\left( L_{n}\geq -r\right) }{\mathbf{P}\left( L_{p}\geq -r\right) }
\end{eqnarray*}%
uniformly in $0\leq x\leq ze^{\varepsilon }\leq N$. Hence we conclude that
given  condition (\ref{Msmalln}) we have as $n\to\infty$%
\begin{eqnarray}
&&\mathbf{P}\left( S_{p}\leq Nc_{p};L_{p,n}\geq xc_{p};L_{n}\geq -r\right)
\notag \\
&&\quad \sim C_{0}\mathbf{P}\left( L_{n}\geq -r\right) \int_{x}^{N}\left(
z-x\right) ^{\alpha \left( 1-\rho \right) }\mathbf{P}\left( S_{p}\in
c_{p}dz|L_{p}\geq -r\right)  \notag \\
&&\quad \sim C_{0}\mathbf{P}\left( L_{n}\geq -r\right) \int_{x}^{N}\left(
z-x\right) ^{\alpha \left( 1-\rho \right) }\mathbb{P}^{(m)}\left( B_{1}\in
dz\right)  \notag \\
&&\quad \qquad =\mathbf{P}\left( L_{n}\geq -r\right) \int_{x}^{N}\left( 1-%
\frac{x}{z}\right) ^{\alpha \left( 1-\rho \right) }\mathbb{P}^{+}\left(
B_{1}\in dz\right) .  \label{ExLim}
\end{eqnarray}%
Using (\ref{SSmall}) and (\ref{ExLim}) to evaluate (\ref{Decomp}) and
letting $N$ to infinity we get (\ref{MinCon}).

For convenience we introduce the notation%
\begin{equation*}
\tau _{n}:=\min \left\{ j:S_{j}=L_{n}\right\} ,\quad \mathcal{A}%
_{u.s.}:=\left\{ Z_{n}>0\text{ for all }n\geq 0\right\}
\end{equation*}%
and recall that by Corollary 1.2 in \cite{AGKV05}, (\ref{tails}) and~(\ref%
{AsH})
\begin{equation}
\mathbf{P}\left( Z_{n}>0\right) \sim \theta \mathbf{P}\left( L_{n}\geq
0\right) \sim \theta l(n)n^{\rho -1}\sim \frac{\theta C_{0}}{V(c_{n})}
\label{DonRen}
\end{equation}%
as $n\rightarrow \infty $, where
\begin{equation}
\theta :=\sum_{k=0}^{\infty }\mathbf{E}[\mathbf{P}_{Z_{k}}^{+}\left(
\mathcal{A}_{u.s.}\right) ;\tau _{k}=k].  \label{DefTheta}
\end{equation}%
Let%
\begin{equation*}
\hat{L}_{k,n}:=\min_{0\leq j\leq n-k}\left( S_{k+j}-S_{k}\right)
\end{equation*}%
and let $\mathcal{\tilde{F}}_{k}$ be the $\sigma -$algebra generated by the
tuple $\left\{ Z_{0},Z_{1},...,Z_{k};Q_{1},Q_{2},...,Q_{k}\right\} $ (see (%
\ref{DefEnvir})).

For further references we formulate two statements borrowed from \cite%
{AGKV05}.

\begin{lemma}
\label{L_cond}(see Lemma 2.5 in \cite{AGKV05}) Assume Condition $A1$. Let $%
Y_{1},Y_{2},...$be a uniformly bounded sequence of random variables adapted
to the filtration $\mathcal{\tilde{F}=}\left\{ \mathcal{\tilde{F}}_{k},k\in
\mathbb{N}\right\} $, which converges $\mathbf{P}^{+}$-a.s. to some random
variable $Y_{\infty }$. Then, as $n\rightarrow \infty $%
\begin{equation*}
\mathbf{E}\left[ Y_{n}|L_{n}\geq 0\right] \rightarrow \mathbf{E}^{+}\left[
Y_{\infty }\right] .
\end{equation*}
\end{lemma}

\begin{lemma}
\label{generaltheorem}(see Lemma 4.1 in \cite{AGKV05}) Assume Condition A1
and let $l\in \mathbb{N}_{0}$. Suppose that $\zeta _{1},\zeta _{2},...$ is a
uniformly bounded sequence of real-valued random variables, which, for every
$k\geq 0$ meets the equality%
\begin{equation*}
\mathbf{E}\left[ \zeta _{n};Z_{k+l}>0,\hat{L}_{k,n}\geq 0\ |\mathcal{\tilde{F%
}}_{k}\right] =\mathbf{P}\left( L_{n}\geq 0\right) \left( \zeta _{k,\infty
}+o(1)\right) \qquad \mathbf{P}\text{-a.s.}
\end{equation*}%
as $n\rightarrow \infty $ with random variables $\zeta _{1,\infty }=\zeta
_{1,\infty }\left( l\right) ,\zeta _{2,\infty }=\zeta _{2,\infty }\left(
l\right) ,....$ Then%
\begin{equation*}
\mathbf{E}\left[ \zeta _{n};Z_{\tau _{n}+l}>0\right] =\mathbf{P}\left(
L_{n}\geq 0\right) \left( \sum_{k=0}^{\infty }\mathbf{E}\left[ \zeta
_{k,\infty };\tau _{k}=k\right] +o(1)\right)
\end{equation*}%
as $n\rightarrow \infty $, where the right-hand side series is absolutely
convergent.
\end{lemma}

For $q\leq p\leq n$ and $u>0$ denote
\begin{equation*}
m(u;n)=m(u;n,p,q):=\min \{q+[u(p-q)],n\}
\end{equation*}%
and set%
\begin{equation*}
\mathcal{X}^{q,p}:=\left\{ X_{u}^{q,p}=e^{-S_{m(u;n)}}Z_{m(u;n)},0\leq
u<\infty \right\} .
\end{equation*}

The next statement is an evident corollary of Theorem 1.3 in \cite{AGKV05}.

\begin{corollary}
\label{C_onedim} Assume Condition $A$. Let $\left( q_{1},p_{1}\right)
,\left( q_{2},p_{2}\right) ,...$ be a sequence of pairs of positive integers
such that $q_{n}\ll p_{n}\ll n$ and $q_{n}\rightarrow \infty $ as $%
n\rightarrow \infty $. Then
\begin{equation*}
\mathcal{L}\left( \mathcal{X}^{q_{n},p_{n}}\ |Z_{n}>0,Z_{0}=1\right)
\Longrightarrow \mathcal{L}\left(W_{u},\ 0\leq u<\infty\right),
\end{equation*}%
where
\begin{equation*}
\mathbf{P}\left( W_{u}=W,\ 0\leq u<\infty\right) =1
\end{equation*}%
for some random variable $W$ such that%
\begin{equation*}
\mathbf{P}\left( 0<W<\infty \right) =1.
\end{equation*}
\end{corollary}

We conclude this section by recalling asymptotic properties of the
distribution of the number of particles in a critical branching process in
random environment at moment $p\ll n$ given $Z_{n}>0.$

\begin{theorem}
\label{T_small} (see \cite{VD2016}) If Condition $A$ is valid, then, as $%
n\rightarrow \infty $
\begin{equation*}
\lim_{n\gg p\rightarrow \infty }\mathbf{P}\left( \frac{\log Z_{p}}{c_{p}}%
\leq z\Big|Z_{n}>0,Z_{0}=1\right) =\mathbb{P}^{+}\left( B_{1}\leq z\right)
\end{equation*}%
for any $z>0$.
\end{theorem}

\section{Reduced processes}

The proof of Theorem \ref{T_reduced} will be divided into several steps
which we formulate as lemmas.

For $f_{n}(s),n=1,2,...,$ specified by (\ref{genfuu}) set

\begin{equation*}
f_{p,n}(s):=f_{p+1}(f_{p+2}(\ldots (f_{n}(s))\ldots )),\;0\leq p\leq
n-1,\;f_{n,n}(s)\equiv 1.
\end{equation*}

We label $Z_{p}$ particles of the $p$th generation by positive numbers $%
1,2,...,Z_{p}$ in an arbitrary but fixed way and denote by $%
Z_{n}^{(i)}(p),\,i=1,2,\ldots ,Z_{p},\,0\leq p\leq n,\,$ the offspring size
at moment $n$ of the population generated by the $i$th particle of the $p$th
generation.

For fixed positive $x$ and $N$ introduce the events
\begin{equation*}
A_{p,n}(x):=\{\,\ln (e+Z_{p,n})\geq xc_{p}\,\},
\end{equation*}

\begin{equation*}
B_{p,n}(N):=\left\{ \left\vert \sum_{i=1}^{Z_{p}}\left( I\left\{
Z_{n}^{(i)}(p)>0\right\} -(1-f_{p,n}(0))\right) \right\vert >\sqrt{%
NZ_{p}(1-f_{p,n}(0))}\right\} ,
\end{equation*}%
and
\begin{equation*}
C_{p,n}:=\left\{ Z_{p}(1-f_{p,n}(0)<e^{\sqrt{c_{p}}}\right\}
\end{equation*}%
and use the notation $\bar{C}_{p,n}$ for the event complementary to $C_{p,n}$%
.

Finally, we use for brevity the notation $\mathbf{P}_{j}(\bullet ):=\mathbf{P%
}(\bullet |Z_{0}=j)$ with the natural agreement that $\mathbf{P}(\bullet ):=%
\mathbf{P}(\bullet |Z_{0}=1)$. In particular, $\mathbf{P}\left(
Z_{n}>0\right) =\mathbf{P}\left( Z_{n}>0|Z_{0}=1\right) $.

\begin{lemma}
\label{L_submain}If the conditions of Theorem \ \ref{T_reduced} are valid
then, for any $j\in \mathbb{N}_{+}$
\begin{equation}
\lim_{N\rightarrow \infty }\lim_{n\gg p\rightarrow \infty }\frac{\mathbf{P}%
_{j}(A_{p,n}(x)B_{p,n}(N);L_{n}\geq 0)}{\mathbf{P}\left( Z_{n}>0\right) }=0
\label{B}
\end{equation}%
and
\begin{equation}
\lim_{n\gg p\rightarrow \infty }\frac{\mathbf{P}_{j}(A_{p,n}(x)C_{p,n};L_{n}%
\geq 0)}{\mathbf{P}\left( Z_{n}>0\right) }=0.  \label{C}
\end{equation}
\end{lemma}

\textbf{Proof.} First we establish the validity of (\ref{B}). To this aim we
temporary introduce the notation $\mathbf{P}^{(\mathcal{F})}(\bullet )$, $%
\mathbf{E}^{(\mathcal{F})}\left[ \bullet \right] $ and $\mathbf{D}^{(%
\mathcal{F})}\left[ \bullet \right] $ for the probability, expectation and
variance calculated for the fixed $\sigma $-algebra $\mathcal{F}=\mathcal{F}%
_{p,n}$ generated by the random variables $\left\{
Z_{0},Z_{1},...,Z_{p}\right\} $ and random probability generating functions $%
f_{1}(s),f_{2}(s),\ldots ,f_{n}(s)$.

First we note that, for $1\leq i\leq Z_{p}$
\begin{eqnarray}
\mathbf{E}^{(\mathcal{F})}\left[ I\left\{ Z_{n}^{(i)}(p)>0\right\} \right]
&=&1-f_{p,n}(0),\quad  \notag \\
\mathbf{D}^{(\mathcal{F})}\left[ I\left\{ Z_{n}^{(i)}(p)>0\right\} \right]
&=&(1-f_{p,n}(0))f_{p,n}(0).  \label{matoz}
\end{eqnarray}%
Besides,
\begin{equation*}
Z_{p,n}=\sum_{i=1}^{Z_{p}}I\left\{ Z_{n}^{(i)}(p)>0\right\} .
\end{equation*}%
Using these relations and applying Chebyshev's inequality to evaluate the
probability under the expectation sign we obtain for sufficiently large $n$
\begin{eqnarray*}
&&\mathbf{P}_{j}(A_{p,n}(x)B_{p,n}(N);L_{n}\geq 0)\leq \mathbf{P}%
_{j}(Z_{p}>0,B_{p,n}(N);L_{n}\geq 0) \\
&&\qquad \qquad =\mathbf{E}\left[ \mathbf{P}_{j}^{(\mathcal{F}%
)}(B_{p,n});Z_{p}>0,L_{n}\geq 0\right] \\
&&\qquad \qquad \leq \mathbf{E}_{j}\left[ \frac{I\left\{ Z_{p}>0,L_{n}\geq
0\right\} }{NZ_{p}(1-f_{p,n}(0)}\mathbf{D}^{(\mathcal{F})}\,\left[
\sum_{i=1}^{Zp}I\left\{ Z_{n}^{(i)}(p)>0\right\} \right] \right] \\
&&\qquad \qquad \qquad \leq N^{-1}\mathbf{P}(L_{n}\geq 0).
\end{eqnarray*}%
This estimate along with (\ref{DonRen}) implies (\ref{B}).

To establish (\ref{C}) observe that
\begin{eqnarray}
&&\mathbf{P}_{j}(A_{p,n}(x)C_{p,n};L_{n}\geq 0)=\mathbf{E}\left[ \mathbf{P}%
_{j}^{(\mathcal{F})}(A_{p,n}(x)C_{p,n});L_{n}\geq 0\right]  \notag \\
&&\qquad \qquad \qquad \leq \frac{1}{x^{2}c_{p}^{2}}\mathbf{E}_{j}\left[
I\left\{ C_{p,n};L_{n}\geq 0\right\} \mathbf{E}^{(\mathcal{F})}[\,\ln
^{2}(e+Z_{p,n})\,]\right] .  \label{estC1}
\end{eqnarray}

Since
\begin{equation*}
(\ln ^{2}(e+x))^{\prime \prime }=\frac{2}{(e+x)^{2}}(1-\ln (e+x)),
\end{equation*}%
the function $\ln ^{2}(e+x)$ is concave on the set $x>0.$ This fact allows
us to apply Jensen's inequality to the internal expectation in the
right-hand side of (\ref{estC1}) and to obtain the estimate
\begin{equation*}
\mathbf{E}^{(\mathcal{F})}\left[ \ln ^{2}(e+Z_{p,n})\right] \leq \ln ^{2}(e+%
\mathbf{E}^{(\mathcal{F})}\left[ Z_{p,n}\right] )=\ln
^{2}(e+Z_{p}(1-f_{p,n}(0))).
\end{equation*}

By the first equality in (\ref{matoz}) we find that, for all sufficiently
large $n,$
\begin{eqnarray*}
&&\mathbf{P}_{j}(A_{p,n}(x)C_{p,n};L_{n}\geq 0)\leq \frac{1}{x^{2}c_{p}^{2}}%
\mathbf{E}_{j}\left[ I\left\{ C_{p,n};L_{n}\geq 0\right\} \ln
^{2}(e+Z_{p}(1-f_{p,n}(0)))\right] \\
&&\qquad \qquad \qquad \qquad \qquad \leq \frac{1}{x^{2}c_{p}^{2}}\ln
^{2}\left( e+e^{\sqrt{c_{p}}}\right) \mathbf{P}\left( L_{n}\geq 0\right)
\leq \frac{K_{1}}{x^{2}c_{p}}\mathbf{P}\left( L_{n}\geq 0\right) .
\end{eqnarray*}

These estimates imply (\ref{C}).

The lemma is proved.

\begin{lemma}
\label{L_joint}If Condition $A$ is valid, then, for any $j=1,2,...$ and $%
x\geq 0$%
\begin{equation*}
\lim_{n\gg p\rightarrow \infty }\mathbf{P}_{j}\left( L_{p,n}\geq
xc_{p}|L_{n}\geq 0,Z_{n}>0\right) =D(x).
\end{equation*}
\end{lemma}

\textbf{Proof}. Clearly,%
\begin{equation}
\mathbf{P}_{j}\left( L_{p,n}\geq xc_{p};L_{n}\geq 0,Z_{n}>0\right) =\mathbf{E%
}\left[ I\left( p,n;x\right) \left( 1-f_{0,n}^{j}(0)\right) \right] ,
\label{New1}
\end{equation}%
where%
\begin{equation*}
I\left( p,n;x\right) :=I\left\{ L_{p,n}\geq xc_{p};L_{n}\geq 0\right\} .
\end{equation*}%
We now select $\gamma >1$ and $l<p$ and write the right-hand side of (\ref%
{New1}) as follows:%
\begin{eqnarray*}
\mathbf{E}\left[ I\left( p,n;x\right) \left( 1-f_{0,n}^{j}(0)\right) \right]
&=&G_{1}\left( m,p,n;x,\gamma \right) +G_{2}\left( p,n;x,\gamma \right)  \\
&&+\,G_{3}\left( m,p,n;x,\gamma \right) ,
\end{eqnarray*}%
where%
\begin{eqnarray*}
G_{1}\left( l,p,n;x,\gamma \right)  &:=&\mathbf{E}\left[ I\left( p,n;x\right)
(f_{0,n}^{j}(0)-f_{0,l}^{j}(0))I\left\{ L_{n\gamma }\geq 0\right\} \right] ,
\\
G_{2}\left( p,n;x,\gamma \right)  &:=&\mathbf{E}\left[ I\left( p,n;x\right)
(1-f_{0,n}^{j}(0))\left( I\left\{ L_{n}\geq 0\right\} -I\left\{ L_{n\gamma
}\geq 0\right\} \right) \right] , \\
G_{3}\left( l,p,n;x,\gamma \right)  &:=&\mathbf{E}\left[ I\left( p,n;x\right)
(1-f_{0,l}^{j}(0))I\left\{ L_{n\gamma }\geq 0\right\} \right] .
\end{eqnarray*}%
By (\ref{IntermedEst}) we have
\begin{eqnarray}
G_{1}\left( l,p,n;x,\gamma \right)  &\leq &\mathbf{E}\left[
(f_{0,n}^{j}(0)-f_{0,l}^{j}(0))I\left\{ L_{n\gamma }\geq 0\right\} \right]
\notag \\
&=&\mathbf{E}\left[ (f_{0,n}^{j}(0)-f_{0,l}^{j}(0))I\left\{ L_{n}\geq
0\right\} \mathbf{P}\left( L_{n(\gamma -1)}^{\ast }\geq -S_{n}|S_{n}\right) %
\right]   \notag \\
&\leq &K_{1}\mathbf{P}\left( L_{n(\gamma -1)}\geq 0\right) \mathbf{E}\left[
(f_{0,n}^{j}(0)-f_{0,l}^{j}(0))I\left\{ L_{n}\geq 0\right\} V\left(
S_{n}\right) \right]   \notag \\
&=&K_{1}\mathbf{P}\left( L_{n(\gamma -1)}\geq 0\right) \mathbf{E}^{+}\left[
f_{0,n}^{j}(0)-f_{0,l}^{j}(0)\right]   \notag \\
&\leq &K_{2}\left( \gamma -1\right) ^{\rho -1}\mathbf{P}\left( L_{n}\geq
0\right) \mathbf{E}^{+}\left[ f_{0,n}^{j}(0)-f_{0,l}^{j}(0)\right] ,
\label{G11}
\end{eqnarray}%
where we have used (\ref{UnifPrecise}) to justify the last inequality. Since
$f_{0,t}(0)\rightarrow f_{0,\infty }(0)\in \left( 0,1\right) $ $\mathbf{P}%
^{+}$ - a.s. as $t\rightarrow \infty ,$ letting in (\ref{G11}) first $%
n\rightarrow \infty $ and then $l\rightarrow \infty ,$ we get in account of (%
\ref{DonRen})%
\begin{equation*}
\lim_{n\rightarrow \infty }\frac{G_{1}\left( l,p,n;x,\gamma \right) }{%
\mathbf{P}\left( L_{n}\geq 0\right) }=0.
\end{equation*}%
Further,
\begin{equation*}
G_{2}\left( p,n;x,\gamma \right) \leq \mathbf{P}\left( L_{n}\geq 0\right) -%
\mathbf{P}\left( L_{n\gamma }\geq 0\right) \leq K_{1}\left( 1-\gamma
^{-(1-\rho )}\right) \mathbf{P}\left( L_{n}\geq 0\right)
\end{equation*}%
implying%
\begin{equation*}
\lim_{\gamma \downarrow 1}\lim_{n\rightarrow \infty }\frac{G_{2}\left(
p,n;x,\gamma \right) }{\mathbf{P}\left( L_{n}\geq 0\right) }=0.
\end{equation*}%
By (\ref{IntermedEst}) we conclude that, given $p\ll n$%
\begin{eqnarray*}
&&G_{3}\left( l,p,n;x,\gamma \right) =\mathbf{E}\left[ (1-f_{0,l}^{j}(0))I%
\left\{ L_{p}\geq 0\right\} \mathbf{P}\left( L_{\gamma n-p}^{\ast }\geq
xc_{p}-S_{p}|S_{p}\right) \right]  \\
&&\qquad \qquad \sim \mathbf{P}\left( L_{n\gamma }\geq 0\right) \mathbf{E}%
\left[ (1-f_{0,l}^{j}(0))I\left\{ L_{p}\geq 0\right\} V\left(
S_{p}-xc_{p}\right) \right]
\end{eqnarray*}%
as $n\rightarrow \infty $. Using Lemma \ref{Renew2} and properties of
regularly varying functions we get, as $p\rightarrow \infty $%
\begin{eqnarray*}
&&\mathbf{E}\left[ (1-f_{0,l}^{j}(0))I\left\{ L_{p}\geq 0\right\} V\left(
S_{p}-xc_{p}\right) \right]  \\
&=&\mathbf{E}\left[ (1-f_{0,l}^{j}(0))I\left\{ L_{p}\geq 0\right\} \frac{%
V\left( S_{p}-xc_{p}\right) }{V\left( c_{p}\right) }\times V\left(
c_{p}\right) I\left\{ \frac{S_{p}}{c_{p}}\geq x\right\} \right]  \\
&\sim &\mathbf{E}\left[ (1-f_{0,l}^{j}(0))I\left\{ L_{p}\geq 0\right\}
\left( \frac{S_{p}}{c_{p}}-x\right) ^{\alpha \left( 1-\rho \right) }V\left(
c_{p}\right) I\left\{ \frac{S_{p}}{c_{p}}\geq x\right\} \right] .
\end{eqnarray*}%
Further,
\begin{eqnarray*}
&&\mathbf{E}\left[ (1-f_{0,l}^{j}(0))I\left\{ L_{p}\geq 0\right\} \left(
\frac{S_{p}}{c_{p}}-x\right) ^{\alpha \left( 1-\rho \right) }V\left(
c_{p}\right) I\left\{ \frac{S_{p}}{c_{p}}\geq x\right\} \right]  \\
&=&\mathbf{E}\left[ (1-f_{0,l}^{j}(0))I\left\{ L_{l}\geq 0\right\} V\left(
c_{p}\right) \mathbf{P}\left( L_{p-l}^{\ast }\geq -S_{l}|S_{l}\right)
\right. \times  \\
&&\times \left. \mathbf{E}\left[ \left( \frac{S_{p-l}^{\ast }+S_{l}}{c_{p}}%
-x\right) ^{\alpha \left( 1-\rho \right) }I\left\{ \frac{S_{p-l}^{\ast
}+S_{l}}{c_{p}}\geq x\right\} |L_{p-l}^{\ast }\geq -S_{l}\right] \right] .
\end{eqnarray*}%
It is not difficult to conclude by Lemma \ref{L_uniform} and (\ref{defD})
that, for any fixed $l$%
\begin{equation*}
\lim_{p\rightarrow \infty }\mathbf{E}\left[ \left( \frac{S_{p-l}^{\ast
}+S_{l}}{c_{p}}-x\right) ^{\alpha \left( 1-\rho \right) }I\left\{ \frac{%
S_{p-l}^{\ast }+S_{l}}{c_{p}}\geq x\right\} |L_{p-l}^{\ast }\geq -S_{l}%
\right] =C_{0}^{-1}D(x).
\end{equation*}%
Besides, for $0\leq x\ll c_{p}$%
\begin{equation*}
V\left( c_{p}\right) \mathbf{P}\left( L_{p-l}^{\ast }\geq -x\right) \sim
V\left( c_{p}\right) \mathbf{P}\left( L_{p-l}^{\ast }\geq 0\right) V(x)\sim
C_{0}V(x)
\end{equation*}%
as $p\rightarrow \infty ,$ leading after evident transformations to
\begin{eqnarray*}
&&\mathbf{E}\left[ (1-f_{0,l}^{j}(0))I\left\{ L_{p}\geq 0\right\} V\left(
S_{p}-xc_{p}\right) \right]  \\
&&\quad \sim D(x)\mathbf{E}\left[ (1-f_{0,l}^{j}(0))I\left\{ L_{l}\geq
0\right\} V(S_{l})\right] =\mathbf{E}^{+}\left[ 1-f_{0,l}^{j}(0)\right] D(x).
\end{eqnarray*}%
Hence we obtain%
\begin{equation*}
\lim_{\gamma \downarrow 1}\lim_{l\rightarrow \infty }\lim_{n\gg p\rightarrow
\infty }\frac{G_{3}\left( l,p,n;x,\gamma \right) }{\mathbf{P}\left(
L_{n}\geq 0\right) }=\mathbf{E}^{+}\left[ 1-f_{0,\infty }^{j}(0)\right] D(x).
\end{equation*}%
To complete the proof of the lemma it remains to note that%
\begin{equation}
\mathbf{P}\left( L_{n}\geq 0,Z_{n}>0,Z_{0}=j\right) =\mathbf{E}\left[ \left(
1-f_{0,n}^{j}(0)\right) |L_{n}\geq 0\right] \mathbf{P}\left( L_{n}\geq
0\right)   \label{Comare}
\end{equation}%
and that%
\begin{equation*}
\lim_{n\rightarrow \infty }\mathbf{E}\left[ \left( 1-f_{0,n}^{j}(0)\right)
|L_{n}\geq 0\right] =\mathbf{E}^{+}\left[ 1-f_{0,\infty }^{j}(0)\right]
\end{equation*}%
according to Lemma \ref{L_cond}.

The lemma is proved.

Set
\begin{equation*}
\eta _{l}:=\frac{f_{l}^{{\prime }{\prime }}(1)}{(f_{l}^{\prime }(1))^{2}}%
,\,l\in \mathbb{N}_{+},
\end{equation*}%
and let%
\begin{eqnarray*}
J^{+}\left( p,r\right)  &:&=\sum_{l=p}^{r-1}\eta
_{l}e^{S_{p}-S_{l}}+e^{S_{p}-S_{r}},\quad J^{-}\left( p,r\right)
:=\sum_{l=p}^{r-1}\eta _{l}e^{S_{r}-S_{l}}, \\
\ \hat{J}^{-}\left( 0,r\right)  &:&=\sum_{l=0}^{r-1}\eta _{l+1}e^{S_{l+1}}.
\end{eqnarray*}%
It is know (see, for instance, Lemma 2.7 in \cite{AGKV05}) that if
Conditions $A1$ and $A2$ are valid then%
\begin{equation*}
J^{+}\left( 0,\infty \right) <\infty \text{ \ }\mathbf{P}^{+}\text{ - a.s.}
\end{equation*}%
and, for any $y>0$
\begin{equation}
\lim_{n\rightarrow \infty }\mathbf{P}\left( J^{+}\left( 0,n\right)
<y|L_{n}\geq 0\right) =\mathbf{P}^{+}\left( J^{+}\left( 0,\infty \right)
<y\right) .  \label{GlobalPlus}
\end{equation}%
In addition (compare with Lemma 2.7 in \cite{AGKV05} or see\ Lemma 6 in \cite%
{VD2005}), if Conditions $A1$ and $A2$ are valid then, for $\mathbf{P}^{-}$
defined in (\ref{DefPminus})
\begin{equation*}
\hat{J}^{-}\left( 0,\infty \right) =\sum_{l=0}^{\infty }\eta
_{l+1}e^{S_{l+1}}<\infty \text{ \ }\mathbf{P}^{-}\text{ - a.s.,}
\end{equation*}%
and for any $y\geq 0$
\begin{equation}
\lim_{n\rightarrow \infty }\mathbf{P}\left( \hat{J}^{-}\left( 0,n\right)
>y|M_{n}\leq 0\right) =\mathbf{P}^{-}\left( \hat{J}^{-}\left( 0,\infty
\right) >y\right) .  \label{GlobalMinus}
\end{equation}

Set%
\begin{equation*}
\tau _{p,n}:=\min \left\{ p\leq j\leq n:S_{j}-S_{p}=\hat{L}_{p,n}\right\} .
\end{equation*}

\begin{lemma}
\label{L_modif}If Condition $A$ is valid, then, for any $j=1,2,...$
\begin{equation}
\lim_{y\rightarrow \infty }\lim_{n\gg p\rightarrow \infty }\mathbf{P}%
_{j}\left( J^{+}\left( \tau _{p,n},n\right) >y|L_{n}\geq 0,Z_{n}>0\right) =0,
\label{Jplus}
\end{equation}%
\begin{equation}
\lim_{y\rightarrow \infty }\lim_{n\gg p\rightarrow \infty }\mathbf{P}%
_{j}\left( J^{-}\left( p,\tau _{p,n}\right) >y|L_{n}\geq 0,Z_{n}>0\right) =0.
\label{Jminus}
\end{equation}
\end{lemma}

\textbf{Proof}. We write%
\begin{eqnarray*}
&&\mathbf{P}_{j}\left( J^{+}\left( \tau _{p,n},n\right) >y;L_{n}\geq
0,Z_{n}>0\right) \\
&&=\mathbf{E}\left[ I\left\{ J^{+}\left( \tau _{p,n},n\right) >y\right\}
\left( 1-f_{0,n}^{j}(0)\right) ;L_{n}\geq 0\right] \\
&&\leq\mathbf{E}\left[ I\left\{ J^{+}\left( \tau _{p,n},n\right) >y\right\}
\left( 1-f_{0,p}^{j}(0)\right) ;L_{n}\geq 0\right] \\
&&=\int_{0}^{\infty }\mathbf{E}\left[ \left( 1-f_{0,p}^{j}(0)\right)
;S_{p}\in c_{p}dz,L_{p}\geq 0\right] \mathbf{P}\left( J^{+}\left( \tau
_{n-p},n-p\right) >y;L_{n-p}\geq -c_{p}z\right) .
\end{eqnarray*}%
Note that
\begin{eqnarray*}
&&\mathbf{P}\left( J^{+}\left( \tau _{n-p},n-p\right) >y;L_{n-p}\geq
-c_{p}z\right) \\
&&\qquad =\sum_{k=0}^{n-p}\mathbf{P}\left( J^{+}\left( k,n-p\right) >y;\tau
_{n-p}=k,L_{n-p}\geq -c_{p}z\right) \\
&&\qquad =\sum_{k=0}^{n-p}\mathbf{P}\left( M_{k}<0;S_{k}\geq -c_{p}z\right)
\mathbf{P}\left( J^{+}\left( 0,n-p-k\right) >y;L_{n-p-k}\geq 0\right) .
\end{eqnarray*}%
In view of (\ref{GlobalPlus}) for any $\varepsilon >0$ there exists $y_{0}$
and $N=N(y_{0},\varepsilon )$ such that, for all $y\geq y_{0}$ and $%
n-p-k\geq N$%
\begin{eqnarray*}
\mathbf{P}\left( J^{+}\left( 0,n-p-k\right) >y,L_{n-p-k}\geq 0\right) &\leq &%
\mathbf{P}\left( J^{+}\left( 0,n-p-k\right) >y_{0},L_{n-p-k}\geq 0\right) \\
&\leq &\mathbf{P}^{+}\left( J^{+}\left( 0,\infty \right) >y_{0}\right)
\mathbf{P}\left( L_{n-p-k}\geq 0\right) \\
&\leq &\varepsilon \mathbf{P}\left( L_{n-p-k}\geq 0\right) .
\end{eqnarray*}%
On the other hand, for each fixed $N$ one can find a sufficiently large $%
y_{1}\geq y_{0}$ such that
\begin{equation*}
\mathbf{P}\left( J^{+}\left( 0,j\right) >y,L_{j}\geq 0\right) \leq
\varepsilon \mathbf{P}\left( L_{j}\geq 0\right)
\end{equation*}%
for all $y\geq y_{1}$. These estimates and (\ref{Comare}) imply%
\begin{eqnarray*}
&&\mathbf{P}\left( J^{+}\left( \tau _{n-p},n-p\right) >y;L_{n-p}\geq
-c_{p}z\right) \\
&&\quad \leq \varepsilon \sum_{k=0}^{n-p}\mathbf{P}\left( M_{k}<0;S_{k}\geq
-c_{p}z\right) \mathbf{P}\left( L_{n-p-k}\geq 0\right) =\varepsilon \mathbf{P%
}\left( L_{n-p}\geq -c_{p}z\right) .
\end{eqnarray*}%
Thus,
\begin{eqnarray*}
&&\int_{0}^{\infty }\mathbf{E}\left[ \left( 1-f_{0,p}^{j}(0)\right)
;S_{p}\in c_{p}dz,L_{p}\geq 0\right] \mathbf{P}\left( J^{+}\left( \tau
_{n-p},n-p\right) >y;L_{n-p}\geq -c_{p}z\right) \\
&&\quad \leq \varepsilon \int_{0}^{\infty }\mathbf{E}\left[ \left(
1-f_{0,p}^{j}(0)\right) ;S_{p}\in c_{p}dz,L_{p}\geq 0\right] \mathbf{P}%
\left( L_{n-p}\geq -c_{p}z\right) \\
&&\quad =\varepsilon \mathbf{E}\left[ \left( 1-f_{0,p}^{j}(0)\right)
;L_{n}\geq 0\right] \leq \varepsilon \mathbf{P}\left( L_{n}\geq 0\right)
\leq \varepsilon K\mathbf{P}_{j}\left( L_{n}\geq 0,Z_{n}>0\right) .
\end{eqnarray*}%
This proves (\ref{Jplus}), since $\varepsilon >0$ may be chosen arbitrary
small.

To prove (\ref{Jminus}) we write%
\begin{eqnarray*}
&&\mathbf{P}_{j}\left( J^{-}\left( p,\tau _{p,n}\right) >y;L_{n}\geq
0,Z_{n}>0\right)  \\
&\leq &\mathbf{E}\left[ \left( 1-f_{0,p}^{j}(0)\right) ;J^{-}\left( p,\tau
_{p,n}\right) >y;L_{n}\geq 0\right]  \\
&=&\int_{0}^{N}\mathbf{E}\left[ \left( 1-f_{0,p}^{j}(0)\right) ;S_{p}\in
c_{p}dz,L_{p}\geq 0\right] \mathbf{P}\left( J^{-}\left( 0,\tau _{n-p}\right)
>y;L_{n-p}\geq -c_{p}z\right)  \\
&&\qquad +\mathbf{E}\left[ \left( 1-f_{0,p}^{j}(0)\right)
;S_{p}>Nc_{p},J^{-}\left( p,\tau _{p,n}\right) >y;L_{n}\geq 0\right] .
\end{eqnarray*}%
By Lemma \ref{C_funct} with $r=0$ for any $\varepsilon >0$ one can find $%
N_{0}$ such that the inequality
\begin{eqnarray}
\mathbf{E}\left[ \left( 1-f_{0,p}^{j}(0)\right) ;S_{p}>Nc_{p},J^{-}\left(
p,\tau _{p,n}\right) >y;L_{n}\geq 0\right]  &\leq &\mathbf{P}\left(
S_{p}>Nc_{p},L_{n}\geq 0\right)   \notag \\
&\leq &\varepsilon \mathbf{P}\left( L_{n}\geq 0\right)   \label{part0}
\end{eqnarray}%
is valid for all $N\geq N_{0}$. Further,%
\begin{eqnarray*}
&&\mathbf{P}\left( J^{-}\left( 0,\tau _{0,n-p}\right) >y;L_{n-p}\geq
-c_{p}z\right)  \\
&&\qquad =\sum_{k=0}^{n-p-1}\mathbf{P}\left( J^{-}\left( 0,k\right)
>y;L_{n-p}\geq -c_{p}z;\tau _{n-p}=k\right)  \\
&&\qquad =\sum_{k=0}^{n-p-1}\mathbf{P}\left( J^{-}\left( 0,k\right)
>y;L_{k}\geq -c_{p}z;\tau _{k}=k\right) \mathbf{P}\left( L_{n-p-k}\geq
0\right)  \\
&&\qquad =\sum_{k=0}^{n-p-1}\mathbf{P}\left( \hat{J}^{-}\left( 0,k\right)
>y;M_{k}\leq 0,S_{k}\geq -c_{p}z\right) \mathbf{P}\left( L_{n-p-k}\geq
0\right) ,
\end{eqnarray*}%
where at the last step we have used the duality principle for random walks.
In view of (\ref{GlobalMinus}) for any $\varepsilon >0$ there exists $%
y_{0}=y_{0}(\varepsilon )$ such that%
\begin{eqnarray*}
\mathbf{P}\left( \hat{J}^{-}\left( 0,k\right) >y;M_{k}\leq 0,S_{k}\geq
-c_{p}z\right)  &\leq &\mathbf{P}\left( \hat{J}^{-}\left( 0,k\right)
>y;M_{k}\leq 0\right)  \\
&\leq &\varepsilon \mathbf{P}\left( M_{k}\leq 0\right)
\end{eqnarray*}%
for all $y\geq y_{0}$. On the other hand, one can show (compare with a
similar statement in \cite{CCh2006} for random walk conditioned to stay
positive) that there exists a proper distribution $G(\cdot )$ with $G(z)\in
\left( 0,1\right) $ for all $z>0$ such that for any $R>0$
\begin{equation*}
\mathbf{P}\left( S_{pR}\geq -c_{p}z|M_{pR}\leq 0\right) \rightarrow
1-G(zR^{1/\alpha })
\end{equation*}%
as $p\rightarrow \infty $ uniformly in $z\in \left[ 0,N\right] $. This leads
to the following chain of estimates being valid for $y\geq y_{0},$ a large
but fixed $T>0$ and $k\leq pT:$
\begin{eqnarray*}
\mathbf{P}\left( \hat{J}^{-}\left( 0,k\right) >y;M_{k}\leq 0,S_{k}\geq
-c_{p}z\right)  &\leq &\mathbf{P}\left( \hat{J}^{-}\left( 0,k\right)
>y;M_{k}\leq 0\right)  \\
\leq \varepsilon \mathbf{P}\left( M_{k}\leq 0\right)  &\leq &\varepsilon
K_{1}\mathbf{P}\left( M_{k}\leq 0,S_{k}\geq -c_{p}z\right) .
\end{eqnarray*}%
This, in account of (\ref{UnifPrecise}) allows us to proceed with one more chain of estimates
\begin{eqnarray}
&&\sum_{k=0}^{pT}\mathbf{P}\left( J^{-}\left( 0,k\right) >y;L_{k}\geq
-c_{p}z;\tau _{k}=k\right) \mathbf{P}\left( L_{n-p-k}\geq 0\right)   \notag
\\
&&\qquad \leq \varepsilon K_{1}\sum_{k=0}^{pT}\mathbf{P}\left( M_{k}\leq
0,S_{k}\geq -c_{p}z\right) \mathbf{P}\left( L_{n-p-k}\geq 0\right)   \notag
\\
&&\qquad \leq \varepsilon K_{1}\mathbf{P}\left( L_{n-p}\geq -c_{p}z\right)
\leq \varepsilon K_{2}V(c_{p}z)\mathbf{P}\left( L_{n-p}\geq 0\right)   \notag
\\
&&\qquad \leq \varepsilon K_{3}z^{\alpha \left( 1-\rho \right) }V(c_{p})%
\mathbf{P}\left( L_{n}\geq 0\right)   \label{part1}
\end{eqnarray}%
being valid for all $z\in \left[ 0,N\right] $.

To consider the case $k\geq Tp$ we use Theorem 4 of \cite{VW09} according to
which%
\begin{equation*}
\mathbf{P}\left( M_{k}\leq 0,S_{k}\geq -c_{p}z\right) \leq K_{1}\frac{%
V(c_{p}z)zc_{p}}{kc_{k}}
\end{equation*}%
if $zc_{p}\leq \varepsilon c_{k}$. Using this bound we get
\begin{eqnarray*}
&&\sum_{k=pT}^{n-1}\mathbf{P}\left( J^{-}\left( 0,k\right) >y;L_{k}\geq
-c_{p}z;\tau _{k}=k\right) \mathbf{P}\left( L_{n-p-k}\geq 0\right)  \\
&&\qquad \leq K_{1}V(c_{p}z)zc_{p}\sum_{k=pT}^{n-1}\frac{1}{kc_{k}}\mathbf{P}%
\left( L_{n-p-k}\geq 0\right) .
\end{eqnarray*}%
By (\ref{asyma}), (\ref{UnifPrecise}) and properties of regularly varying
functions we conclude that%
\begin{eqnarray*}
\sum_{k=pT}^{n-1}\frac{1}{kc_{k}}\mathbf{P}\left( L_{n-p-k}\geq 0\right)
&=&\sum_{k=pT}^{n/2+p}\frac{1}{kc_{k}}\mathbf{P}\left( L_{n-p-k}\geq
0\right) +\sum_{k=n/2+p+1}^{n-1}\frac{1}{kc_{k}}\mathbf{P}\left(
L_{n-p-k}\geq 0\right)  \\
&\leq &K_{1}\mathbf{P}\left( L_{n}\geq 0\right) \sum_{k=pT}^{n/2+p}\frac{1}{%
kc_{k}}+\frac{K_{2}}{nc_{n}}\sum_{j=0}^{n/2}\mathbf{P}\left( L_{j}\geq
0\right)  \\
&\leq &\frac{K_{3}}{c_{pT}}\mathbf{P}\left( L_{n}\geq 0\right) +\frac{K_{4}}{%
nc_{n}}n\mathbf{P}\left( L_{n}\geq 0\right)  \\
&\leq &\frac{K_{5}}{c_{pT}}\mathbf{P}\left( L_{n}\geq 0\right) .
\end{eqnarray*}%
As a result we get%
\begin{eqnarray}
&&\sum_{k=pT}^{n-1}\mathbf{P}\left( J^{-}\left( 0,k\right) >y;L_{k}\geq
-c_{p}z;\tau _{k}=k\right) \mathbf{P}\left( L_{n-p-k}\geq 0\right)   \notag
\\
&&\quad \leq K_{1}V(c_{p}z)zc_{p}\frac{K_{5}}{c_{pT}}\mathbf{P}\left(
L_{n}\geq 0\right) \leq K_{6}z^{\alpha \left( 1-\rho \right) +1}\frac{c_{p}}{%
c_{pT}}V(c_{p})\mathbf{P}\left( L_{n}\geq 0\right)   \notag \\
&&\quad \leq K_{7}z^{\alpha \left( 1-\rho \right) +1}\frac{1}{T^{1/\alpha }}%
V(c_{p})\mathbf{P}\left( L_{n}\geq 0\right) .  \label{part2}
\end{eqnarray}%
Combining (\ref{part1}) and (\ref{part2}) we see that
\begin{equation*}
\mathbf{P}\left( J^{-}\left( 0,\tau _{0,n-p}\right) >y;L_{n-p}\geq
-c_{p}z\right) \leq K_{8}z^{\alpha \left( 1-\rho \right) }\left( \frac{z}{%
T^{1/\alpha }}+\varepsilon \right) V(c_{p})\mathbf{P}\left( L_{n}\geq
0\right) .
\end{equation*}%
Thus, assuming that $N\leq \varepsilon T^{1/\alpha }$ we get in account of (%
\ref{AsH})
\begin{eqnarray*}
&&\int_{0}^{N}\mathbf{E}\left[ \left( 1-f_{0,p}^{j}(0)\right) ;S_{p}\in
c_{p}dz,L_{p}\geq 0\right] \mathbf{P}\left( J^{-}\left( 0,\tau _{n-p}\right)
>y;L_{n-p}\geq -c_{p}z\right)  \\
&&\quad \leq \varepsilon K_{1}\mathbf{P}\left( L_{n}\geq 0\right)
V(c_{p})\times  \\
&&\quad \times \int_{0}^{N}\mathbf{E}\left[ \left( 1-f_{0,p}^{j}(0)\right)
;S_{p}\in c_{p}dz,L_{p}\geq 0\right] z^{\alpha \left( 1-\rho \right) }\left(
\frac{z}{T^{1/\alpha }}+\varepsilon \right)  \\
&&\quad \leq \varepsilon K_{1}\mathbf{P}\left( L_{n}\geq 0\right)
V(c_{p})\int_{0}^{N}z^{\alpha \left( 1-\rho \right) }\mathbf{P}\left(
S_{p}\in c_{p}dz,L_{p}\geq 0\right)  \\
&&\quad \leq \varepsilon K_{1}\mathbf{P}\left( L_{n}\geq 0\right) V(c_{p})%
\mathbf{P}\left( L_{p}\geq 0\right) \int_{0}^{N}z^{\alpha \left( 1-\rho
\right) }\mathbf{P}\left( S_{p}\in c_{p}dz|L_{p}\geq 0\right)  \\
&&\quad \leq \varepsilon K_{2}C_{0}\mathbf{P}\left( L_{n}\geq 0\right)
\mathbb{P}^{+}\left( B_{1}\leq N\right) \leq \varepsilon K_{3}\mathbf{P}%
\left( L_{n}\geq 0\right) .
\end{eqnarray*}%
This estimate combined with (\ref{part0}) proves (\ref{Jminus}).

\begin{lemma}
\label{L_main}If Condition $A$ is valid, then, for any $j=1,2,...$ and $%
x\geq 0$%
\begin{equation*}
\lim_{n\gg p\rightarrow \infty }\mathbf{P}_{j}\left( A_{p,n}(x)|L_{n}\geq
0,Z_{n}>0\right) =D(x\mathbb{)}.
\end{equation*}
\end{lemma}

\textbf{Proof}. It follows from Lemma \ref{L_submain} and the inequalities%
\begin{eqnarray*}
&&\mathbf{P}_{j}(A_{p,n}(x)\bar{B}_{p,n}(N)\bar{C}_{p,n};L_{n}\geq 0)\leq
\mathbf{P}_{j}(A_{p,n}(x);L_{n}\geq 0) \\
&&\qquad \quad \leq \mathbf{P}_{j}(A_{p,n}(x)\bar{B}_{p,n}(N)\bar{C}%
_{p,n};L_{n}\geq 0,Z_{n}>0) \\
&&\qquad \quad +\mathbf{P}_{j}(A_{p,n}(x)C_{p,n};L_{n}\geq 0)+\mathbf{P}%
_{j}(A_{p,n}(x)B_{p,n}(N);L_{n}\geq 0)
\end{eqnarray*}%
that, in fact, we need to show that
\begin{equation*}
\lim_{N\rightarrow \infty }\lim_{n\gg p\rightarrow \infty }\frac{\mathbf{P}%
_{j}(A_{p,n}(x)\bar{B}_{p,n}(N)\bar{C}_{p,n};L_{n}\geq 0,Z_{n}>0)}{\mathbf{P}%
_{j}(Z_{n}>0,L_{n}\geq 0)}=D(x\mathbb{)}.
\end{equation*}%
Using the equality
\begin{equation*}
Z_{p,n}=Z_{p}(1-f_{p,n}(0))+\sum_{l=1}^{Z_{p}}\left( I\left\{
Z_{n}^{(l)}(p)>0\right\} -(1-f_{p,n}(0))\,\right) ,
\end{equation*}%
the estimate
\begin{equation*}
\left\vert \,\sum_{l=1}^{Z_{p}}\left( I\left\{ Z_{n}^{(l)}(p)>0\right\}
-(1-f_{p,n}(0))\right) \right\vert \leq \sqrt{NZ_{p}(1-f_{p,n}(0)},
\end{equation*}%
being valid on the set $\bar{B}_{p,n}(N),$ and recalling (\ref{B}), (\ref{C}%
) and the fact that
\begin{equation*}
Z_{p}(1-f_{p,n}(0)))\geq e^{\sqrt{c_{p}}}
\end{equation*}%
on the set $\bar{C}_{p,n}\cap \{Z_{p}>0\}$, we conclude that for any $\delta
>0$ there exists a number $n_{0}=n_{0}(\delta )$ such that for $n\geq n_{0}$
\begin{eqnarray*}
&&\mathbf{P}_{j}(A_{p,n}(x)\bar{B}_{p,n}(N)\bar{C}_{p,n},L_{n}\geq 0,Z_{n}>0)
\\
&\leq &\mathbf{P}_{j}\left( \ln (e+Z_{p}(1-f_{p,n}(0))(1+\delta ))\geq
xc_{p};\bar{C}_{p,n};\,L_{n}\geq 0,Z_{n}>0\right) +\alpha _{1}(p,n) \\
&\leq &\mathbf{P}_{j}\left( \ln Z_{p}+\ln (1-f_{p,n}(0))+\ln (1+2\delta
)\geq xc_{p};L_{n}\geq 0,Z_{n}>0\right) +\alpha _{1}(p,n),
\end{eqnarray*}%
where%
\begin{equation*}
\lim_{n\gg p\rightarrow \infty }\frac{\left\vert \alpha _{1}(p,n)\right\vert
}{\mathbf{P}_{j}\left( Z_{n}>0\right) }=0.
\end{equation*}%
In view of
\begin{equation*}
1-f_{p,n}(0)\leq f_{p}^{\prime }(1)\left( 1-f_{p+1,n}(0)\right) \leq
\min_{p\leq j\leq n}\prod\limits_{p\leq i\leq j}f_{i}^{\prime }(1)=e^{\hat{L}%
_{p,n}},
\end{equation*}%
we get%
\begin{eqnarray*}
&&\mathbf{P}_{j}\left( \ln Z_{p}+\ln (1-f_{p,n}(0))+\ln (1+2\delta )\geq
xc_{p};L_{n}\geq 0,Z_{n}>0\right) \\
&&\quad \leq \mathbf{P}_{j}\left( \frac{1}{c_{p}}\ln \frac{Z_{p}}{e^{S_{p}}}+%
\frac{1}{c_{p}}L_{p,n}+\ln (1+2\delta )\geq x;L_{n}\geq 0,Z_{n}>0\right) .
\end{eqnarray*}%
Using the equivalences
\begin{eqnarray*}
\mathbf{P}_{j}\left( L_{n}\geq 0,Z_{n}>0\right) &\sim &\mathbf{P}_{j}\left(
Z_{n}>0|L_{n}\geq 0\right) \mathbf{P}\left( L_{n}\geq 0\right) \\
&\sim &\frac{1}{\theta }\mathbf{E}^{+}\left( 1-f_{0,\infty }^{j}(0)\right)
\mathbf{P}\left( Z_{n}>0\right)
\end{eqnarray*}%
valid as $n\rightarrow \infty $, and following the proof of Theorem 1.1 in
\cite{AGKV05}, we obtain%
\begin{eqnarray*}
&&\mathbf{P}_{j}(A_{p,n}(x)\bar{B}_{p,n}(N)\bar{C}_{p,n}|L_{n}\geq 0,Z_{n}>0)
\\
&&\qquad \qquad \leq \mathbf{P}_{j}\left( \frac{1}{c_{p}}\ln \frac{Z_{p}}{%
e^{S_{p}}}+\frac{1}{c_{p}}L_{p,n}+\ln (1+2\delta )\geq x|L_{n}\geq
0,Z_{n}>0\right) .
\end{eqnarray*}%
This inequality, Corollary \ref{C_onedim} and Lemma \ref{L_joint} yield
\begin{eqnarray*}
&&\lim_{N\rightarrow \infty }\limsup_{n\gg p\rightarrow \infty }\mathbf{P}%
_{j}\left( A_{p,n}(x)\bar{B}_{p,n}(N)\bar{C}_{p,n}~|~L_{n}\geq
0,Z_{n}>0\right) \\
&&\quad \leq \lim_{\delta \downarrow 0}\lim_{n\gg p\rightarrow \infty }%
\mathbf{P}_{j}\left( \frac{1}{c_{p}}L_{p,n}+\ln (1+3\delta )\geq
x~|~L_{n}\geq 0,Z_{n}>0\right) =D(x).
\end{eqnarray*}%
To get a similar estimate from below observe that according to relations
(2.2) and (2.3) of \cite{GK2000}
\begin{eqnarray*}
1-f_{p,n}(0) &\geq &\left( \sum_{l=p}^{n-1}\eta
_{l}e^{-(S_{l}-S_{p})}+e^{-(S_{n}-S_{p})}\right) ^{-1} \\
&=&e^{\hat{L}_{p,n}}\left( \sum_{l=p}^{n-1}\eta _{l}e^{-(S_{l}-\hat{L}%
_{p,n})}+e^{-(S_{n}-\hat{L}_{p,n})}\right) ^{-1} \\
&=&e^{\hat{L}_{p,n}}\left( J^{-}\left( p,\tau _{p,n}\right) +J^{+}\left(
\tau _{p,n},n\right) \right) ^{-1}.
\end{eqnarray*}

Thus,%
\begin{equation*}
\log \left( 1-f_{p,n}(0)\right) \geq \hat{L}_{p,n}-\log \left( J^{-}\left(
p,\tau _{p,n}\right) +J^{+}\left( \tau _{p,n},n\right) \right) .
\end{equation*}%
According to Lemma \ref{L_modif}
\begin{equation*}
\lim_{y\rightarrow \infty }\lim_{n\gg p\rightarrow \infty }\mathbf{P}%
_{j}\left( \log \left( J^{-}\left( p,\tau _{p,n}\right) +J^{+}\left( \tau
_{p,n},n\right) \right) >y|L_{n}\geq 0,Z_{n}>0\right) =0.
\end{equation*}%
Hence for any $\delta >0$ we get for all sufficiently large $p$ and $n$
\begin{eqnarray*}
&&\mathbf{P}_{j}\left( \ln Z_{p}+\ln (1-f_{p,n}(0))+\ln (1-\delta )\geq
xc_{p};L_{n}\geq 0,Z_{n}>0\right) \\
&&\qquad \quad \geq \mathbf{P}_{j}\left( \frac{1}{c_{p}}\ln \frac{Z_{p}}{%
e^{S_{p}}}+\frac{1}{c_{p}}L_{p,n}+\ln (1-2\delta )\geq x;L_{n}\geq
0,Z_{n}>0\right)
\end{eqnarray*}%
leading by Lemma \ref{L_joint} to the following estimate from below:%
\begin{eqnarray*}
&&\lim \inf_{n\gg p\rightarrow \infty }\mathbf{P}_{j}\left( A_{p,n}(x)\bar{B}%
_{p,n}(N)\bar{C}_{p,n}~|~L_{n}\geq 0,Z_{n}>0\right) \\
&\geq &\lim_{\delta \downarrow 0}\lim_{n\gg p\rightarrow \infty }\mathbf{P}%
_{j}\left( \frac{1}{c_{p}}L_{p,n}+\ln (1-3\delta )\geq x~|~L_{n}\geq
0,Z_{n}>0\right) =D(x\mathbb{)}.
\end{eqnarray*}%
Lemma \ref{L_main} is proved.

\textbf{Proof of Theorem \ref{T_reduced}. } First we note that to check the
validity of (\ref{One_dim}) it is sufficent to investigate the asymptotic
behavior of the probability of the event $A_{p,n}(x)$. We use Lemma \ref%
{generaltheorem} to this aim. For $z,p,n\in \mathbb{N}_{0}$ with $p\leq n$
set%
\begin{equation*}
\psi (z,p,n):=\mathbf{P}_{z}\left( A_{p,n}(x),L_{n}\geq 0\right) .
\end{equation*}%
Clearly, $\psi (0,p,n)=0$. We know by Lemmas \ref{L_main} and \ref{L_cond}
that if $n\gg p=p(n)\rightarrow \infty $ then
\begin{eqnarray}
\psi (z,p,n) &\sim &\mathbf{P}_{z}\left( L_{n}\geq 0,Z_{n}>0\right) D(x)
\notag \\
&=&\mathbf{P}\left( L_{n}\geq 0\right) \mathbf{E}\left[
1-f_{0,n}^{z}(0)|L_{n}\geq 0\right] D(x)  \notag \\
&\sim &\mathbf{P}\left( L_{n}\geq 0\right) \mathbf{E}^{+}\left[
1-f_{0,\infty }^{z}(0)\right] D(x)  \notag \\
&=&\mathbf{P}\left( L_{n}\geq 0\right) \mathbf{P}_{z}^{+}\left( \mathcal{A}%
_{u.s}\right) D(x).  \label{Prelim1}
\end{eqnarray}%
In addition, for $k\leq p\leq n$%
\begin{equation}
\mathbf{E}\left[ I\left\{ A_{p,n}(x)\right\} ,L_{k,n}\geq 0|\mathcal{\tilde{F%
}}_{k}\right] =\psi (Z_{k},p-k,n-k).  \label{Prelim2}
\end{equation}%
Relations (\ref{Prelim1}) and (\ref{Prelim2}) show that we may apply Lemma %
\ref{generaltheorem} to
\begin{equation*}
\zeta _{n}:=I\left\{ A_{p,n}(x)\right\} ,\ \zeta _{k,\infty }:=\mathbf{P}%
_{Z_{k}}^{+}\left( \mathcal{A}_{u.s}\right) D(x)
\end{equation*}%
and $l=0$ to conclude that%
\begin{eqnarray*}
\mathbf{P}\left( A_{p,n}(x)\right)  &=&\mathbf{P}\left(
A_{p,n}(x);Z_{n}>0\right)  \\
&\sim &D(x)\theta \mathbf{P}\left( L_{n}\geq 0\right) \sim D(x)\mathbf{P}%
\left( Z_{n}>0\right) .
\end{eqnarray*}%
This completes the proof of (\ref{One_dim})\textbf{. }

\end{document}